\documentclass[a4paper, 11pt]{article}

\usepackage{amsmath}
\usepackage{amssymb}
\usepackage{graphicx}

\usepackage[margin=1in]{geometry}

\newtheorem{theorem}{Theorem}[section]
\newtheorem{lemma}{Lemma}[section]
\newtheorem{proposition}{Proposition}[section]
\newtheorem{corollary}{Corollary}[section]

\begin{document}

\title{Diagonal Riccati Stability and the Hadamard Product}

\author{Alexander Aleksandrov\thanks{Saint Petersburg State University, Saint Petersburg, Russia}, Oliver Mason \thanks{Dept. of Mathematics and Statistics/Hamilton Institute, Maynooth University, Co. Kildare, Ireland \& Lero, The Irish Software Research Centre} \footnote{Corresponding author, email: oliver.mason@nuim.ie} and Anna Vorob'eva\thanks{Saint Petersburg State University, Saint Petersburg, Russia}}

\maketitle

\begin{abstract}
We first present an extension of a recent characterisation of diagonal Riccati stability and, using this, extend a result of Kraaijevanger on diagonal Lyapunov stability to Riccati stability of time-delay systems.  We also describe a class of transformations that preserve the property of being diagonally Riccati stable and apply these two results to provide novel stability results for classes of time-delay systems. 
\end{abstract}


\vspace{5mm}


\section{Background and Introduction}
The problem of Riccati stability was introduced in \cite{VERR} and is motivated by the stability theory of linear time-delay systems.  Formally, a pair $(A, B)$ is said to be Riccati stable if there exist $P = P^T \succ 0$, $Q = Q^T \succ 0$ such that
\begin{equation}
\label{eq:Ricc3}
A^TP+PA + Q + PBQ^{-1}B^TP \prec 0,
\end{equation}
where $M \prec 0$ ($M  \succ  0$) denotes that the matrix $M=M^T$ is negative definite (positive definite).  Throughout the paper, $M \succeq 0$ ($M \preceq 0$) denotes that $M$ is positive semi-definite (negative semi-definite).  When matrices $P$, $Q$ satisfying \eqref{eq:Ricc3} exist they define a quadratic Lyapunov-Krasovskii functional establishing stability for the time-delay system 
\begin{equation}\label{eq:tdelay}
\dot{x}(t) = Ax(t) + B x(t - \tau),
\end{equation}
where $\tau \geq 0$ can be any fixed nonnegative delay.  

$A \in \mathbb{R}^{n \times n}$ is \emph{Metzler} if $a_{ij} \geq 0$ for $i \neq j$.  We denote the spectrum of $A$ by $\sigma(A)$ and the spectral abscissa of $A$ by $\mu(A)$: formally,
$$\mu(A) := \max \{\textrm{Re}(\lambda) \mid \lambda \in \sigma(A)\}$$
and say that $A$ is Hurwitz if $\mu(A) < 0$.  

We denote the standard basis of $\mathbb{R}^n$ by $e_1, \ldots, e_n$ and we use $\mathbf{1}_n$ to denote the vector in $\mathbb{R}^n$, all of whose entries are equal to one.  For $A \in \mathbb{R}^{n \times n}$, $\textrm{diag}(A)$ is the vector $v$ in $\mathbb{R}^n$ with $v_i = a_{ii}$, $1 \leq i \leq n$.  $\textrm{Sym}(n, \mathbb{R})$ denotes the space of $n \times n$ symmetric matrices with real entries.  

For a real number $x$, $\textrm{sign}(x)$ is given by $+1$ if $x \geq 0$ and $-1$ for $x < 0$ respectively.  

For vectors $v, w$ in $\mathbb{R}^n$, we write: $v \geq w$ if $v_i \geq w_i$ for $1 \leq i \leq n$; $v > w$ if $v \geq w$, $v \neq w$; $v \gg w$ if $v_i > w_i$ for $1 \leq i \leq n$. 

When diagonal positive definite solutions $P, Q$ of \eqref{eq:Ricc3} exist, we say that the pair $(A, B)$ is \emph{diagonally Riccati stable}. 

In \cite{AlexMas16}, a necessary and sufficient condition for a given pair $(A, B)$ of matrices in $\mathbb{R}^{n \times n}$ to be diagonally Riccati stable was described.  This result extended naturally the celebrated condition of Barker, Berman and Plemmons for diagonal Lyapunov stability \cite{BBP, Bob1}.  The existence of diagonal solutions to \eqref{eq:Ricc3} allows the construction of Lyapunov-Krasovskii functionals of particularly simple form.  As with the case of diagonal Lyapunov stability for undelayed systems \cite{KasBha}, such functionals prove useful in establishing absolute stability conditions for classes of nonlinear time-delay systems \cite{AlexMas16} (see \cite{AlexMas14} for corresponding work on discrete time systems).

Formally, the following result was proven in \cite{AlexMas16}. 
\begin{theorem}
\label{thm:TM} Let $A, B \in \mathbb{R}^{n \times n}$ be given.  The following are equivalent.
\begin{itemize}
\item[(i)] There exist $P \succ 0$, $Q \succ 0$ diagonal satisfying \eqref{eq:Ricc3}.
\item[(ii)] For every non-zero positive semi-definite
\begin{equation}
\label{eq:H}
H = \left(\begin{array}{c c}
		H_{11} & H_{12} \\
		H_{12}^T & H_{22}
	  \end{array}\right)
\end{equation}
in $Sym(2n, \mathbb{R})$ with $\textrm{diag}(H_{11}) \geq \textrm{diag}(H_{22})$, the matrix
$$AH_{11} + BH_{12}^T$$
has a negative diagonal entry.
\end{itemize}
\end{theorem}

For $H$ of the form \eqref{eq:H}, we use $h^{ab}_{ij}$ to denote the $i, j$ element of $H_{ab}$ for $1 \leq a, b \leq 2$, $1 \leq i, j \leq n$.

Using Theorem \ref{thm:TM}, necessary and sufficient conditions were derived for diagonal Riccati stability for pairs $(A, B)$ where: (i) $A$ is Metzler and $B$ is nonnegative; (ii) $A$ and $B$ are both upper (lower) triangular.  

In particular, it was shown that in case (i), diagonal Riccati stability is equivalent to the condition that $A+B$ is Hurwitz.  This established the existence of a diagonal Lyapunov functional for asymptotically stable positive linear time-delay systems, providing a natural extension of a fundamental property of positive linear time-invariant (LTI) systems \cite{FarRin}.  Furthermore, this fact strengthened the main result of \cite{MAS1} which showed that under the same condition ($A+B$ Hurwitz), there exists a diagonal $P \succ 0$ and $Q \succ 0$ (not necessarily diagonal) satisfying \eqref{eq:Ricc3}.  

In the current paper, we derive a slight extension of Theorem \ref{thm:TM} to identify further classes of matrices for which simple necessary and sufficient conditions for diagonal Riccati stability can be derived.  We will also describe a characterisation of diagonal Riccati stability in terms of the so-called Hadamard product, extending a well known result of \cite{KRA} for diagonal Lyapunov stability.  We will use this alternative characterisation to derive necessary and sufficient conditions for diagonal Riccati stability for two classes of pairs of matrices in $\mathbb{R}^{3 \times 3}$.

\section{Diagonal Riccati Stability and Hadamard Products}
For some of our later results, we will need the following slight extension of Theorem \ref{thm:TM}.
\begin{theorem}
\label{thm:TM2} Let $A, B \in \mathbb{R}^{n \times n}$ be given.  The following are equivalent.
\begin{itemize}
\item[(i)] There exist $P \succ 0$, $Q \succ 0$ diagonal satisfying \eqref{eq:Ricc3}.
\item[(iia)] For every non-zero positive semi-definite $H$ given by \eqref{eq:H}
in $Sym(2n, \mathbb{R})$ with $\textrm{diag}(H_{11}) = \textrm{diag}(H_{22})$, the matrix
$$AH_{11} + BH_{12}^T$$
has a negative diagonal entry.
\end{itemize}
\end{theorem}
\textbf{Proof:} It is enough to show that condition (iia) above is equivalent to condition (ii) in Theorem \ref{thm:TM}.  The implication (ii) $\Rightarrow$ (iia) is trivial.  Assume (iia) and let $H$ satisfy $\textrm{diag}(H_{11}) \geq \textrm{diag}(H_{22})$.  Write $e = \textrm{diag}(H_{11}) - \textrm{diag}(H_{22})$ and set $E$ to be the corresponding diagonal matrix.  Then the matrix 
$$\hat{H} = H + \left(\begin{array}{c c}
							0 & 0\\
                            0 & E
                        \end{array}\right)$$
clearly satisfies (iia) and is also positive semi-definite.  So we can conclude that $$AH_{11} + BH_{12}^T$$ has a negative diagonal entry.  This completes the proof. 

Given two matrices $A \in \mathbb{R}^{n \times n}$, $B \in \mathbb{R}^{n \times n}$, the Hadamard product $A \circ B$ is the matrix $C$ with $(i, j)$ entry given by $a_{ij}b_{ij}$. 

\begin{proposition}
\label{prop:SHad1} Let $(A, B)$ be a diagonally Riccati stable pair and let $$S = \left(\begin{array}{c c} S_{11} & S_{12} \\
						S_{12}^T & S_{22}
                        \end{array}\right)$$
be positive semi-definite, with $\textrm{diag}(S_{11}) = \textrm{diag}(S_{22}) \gg 0$.  Then the pair $(A \circ S_{11}, B \circ S_{12})$ is also diagonally Riccati stable. 
\end{proposition}
\textbf{Proof:} Let $H$ be a non-zero positive semi-definite matrix in $\textrm{Sym}(2n, \mathbb{R})$ given by \eqref{eq:H}.  It follows from the Schur product theorem that the matrix $G =  S \circ H$ is positive semi-definite; as all diagonal entries of $S$ are positive, $G$ is non-zero.  Moreover, if $\textrm{diag}(H_{11}) =\textrm{diag}(H_{22})$ and $\textrm{diag}(S_{11}) =\textrm{diag}(S_{22})$, we also have that for the matrix $G$, $\textrm{diag}(G_{11}) =\textrm{diag}(G_{22})$.  As the pair $(A, B)$ is diagonally Riccati stable, we can conclude that $A G_{11} + B G_{12}^T$ has a negative diagonal entry. 

Clearly, $G_{11} = S_{11} \circ H_{11}$ and $G_{12} = S_{12} \circ H_{12}$.  For any $i$, $1 \leq i \leq n$, 
\begin{eqnarray*}
[A(S_{11}\circ H_{11})]_{ii} &=& \sum_{j=1}^n a_{ij}(S_{11}\circ H_{11})_{ji}\\
&=& \sum_{j=1}^n a_{ij}s^{11}_{ji}h^{11}_{ji}\\
&=& \sum_{j=1}^n a_{ij}s^{11}_{ij}h^{11}_{ji}\\
&=& \sum_{j=1}^n (A\circ S_{11})_{ij} h^{11}_{ji}\\
&=& [(A\circ S_{11})H_{11}]_{ii}.
\end{eqnarray*}
A similar calculation reveals that 
\begin{eqnarray*}
[B(S_{12}^T\circ H_{12}^T)]_{ii} &=& \sum_{j=1}^n b_{ij}(S_{12}^T\circ H_{12}^T)_{ji}\\
&=& \sum_{j=1}^n b_{ij}s^{12}_{ij}h^{12}_{ij}\\
&=& \sum_{j=1}^n (B\circ S_{12})_{ij} h^{12}_{ij}\\
&=& [(B\circ S_{12})H_{12}^T]_{ii}.
\end{eqnarray*}
As $A G_{11} + B G_{12}^T$ has a negative diagonal entry, it follows that for some $i$, 
$$[A(S_{11}\circ H_{11})]_{ii} + [B(S_{12}^T\circ H_{12}^T)]_{ii} < 0$$
and therefore that
$$[(A\circ S_{11})H_{11}]_{ii}+[(B\circ S_{12})H_{12}^T]_{ii} < 0.$$
It follows from Theorem \ref{thm:TM2} that the pair $(A \circ S_{11}, B \circ S_{12})$ is diagonally Riccati stable as claimed. 

\emph{Remark on D-stability and diagonal stability}

The concept of D-stability was originally motivated by applications in areas such as Ecology and Economics and, for the case of a linear time-invariant (LTI) system with system matrix $A$, this property amounts to requiring that $DA$ is Hurwitz for every diagonal $D \succ 0$.  One natural way of extending this to the time-delay case is to require that the time-delay system with matrices $DA, DB$ is stable for all diagonal $D \succ 0$.  It is not a difficult calculation to see that this will hold if the pair $(A, B)$ is diagonally Riccati stable.

The main result of this section provides a direct extension of Theorem 1.2 of \cite{KRA} to the setting of diagonal Riccati stability.  We first recall that a matrix $A$ in $\mathbb{R}^{n \times n}$ is a P matrix if every principal minor of $P$ is positive.  It is a classical result that this is equivalent to the condition that for every non-zero $x$ in $\mathbb{R}^n$, there is some index $i$ with $x_i (Ax)_i > 0$.  

Before stating the main result of this section, we recall a bassic fact relating Lyapunov diagonal stability to the P property for a single matrix \cite{KasBha}. 
\begin{lemma}
\label{lem:Pdiag} Let $A \in \mathbb{R}^{n \times n}$ be diagonally stable.  Then $-A$ is a P-matrix.
\end{lemma}

\begin{theorem}
\label{thm:Kra} Let $A, B$ in $\mathbb{R}^{n \times n}$ be given.  The pair $(A, B)$ is diagonally Riccati stable if and only if $-(A\circ S_{11} + B \circ S_{12})$ is a P-matrix for all 
\begin{equation}
\label{eq:HadS} S = \left(\begin{array}{c c}
							S_{11} & S_{12}\\
                            S_{12}^T & S_{22}
                            \end{array}\right)
\end{equation}
satisfying $S \succeq 0$, $\textrm{diag}(S_{11}) = \textrm{diag}(S_{22}) \gg 0$.  
\end{theorem}
\textbf{Proof:} First suppose that $(A, B)$ is a diagonally Riccati stable pair.  It follows from Proposition \ref{prop:SHad1} that for all $S$ satisfying the hypotheses of the theorem, that $(A\circ S_{11}, B \circ S_{12})$ is also diagonally Riccati stable.  This implies that the matrix 
$$L:=A\circ S_{11} + B \circ S_{12}$$ is diagonally Lyapunov stable and hence $-L$ is a P-matrix.  

For the converse, let $H \succeq 0$ be a non-zero matrix in $Sym(2n, \mathbb{R})$ of the form \eqref{eq:H} with $\textrm{diag}(H_{11}) = \textrm{diag}(H_{22})$.  Define the matrix $S \in Sym(2n, \mathbb{R})$ by setting:
$$ s_{ij} = \begin{cases}
			h_{ij} & \mbox{ if } i \neq j\\
            h_{ii} & \mbox{ if } i=j \& h_{ii}> 0\\
            1 & \mbox{ if } i=j \& h_{ii} = 0.
		\end{cases}
$$
It can then be verified that $S$ is positive semi-definite and that $s_{ii} > 0$ for $1 \leq i \leq 2n$.  Moreover, as $\textrm{diag}(H_{11}) = \textrm{diag}(H_{22})$, it follows that the same holds for the diagonal elements of $S$ so that $S$ satisfies the hypotheses of the theorem.  Thus we can conclude that 
$$-(A\circ S_{11} + B \circ S_{12})$$ is a P-matrix. Now choosing $x \in \mathbb{R}^n$ with $x_i = 1$ if $h_{ii} > 0$ and $x_i = 0$ if $h_{ii} = 0$ for $1\leq i \leq n$, we conclude that there is some index $i$ with 
$$x_i [(A\circ S_{11} + B \circ S_{12})x]_i < 0.$$
If we expand this we find that
$$\sum_{j=1}^n a_{ij}h^{11}_{ij} + \sum_{j=1}^n b_{ij}h^{12}_{ij} < 0$$
or equivalently 
$$(AH_{11} + B H_{12}^T)_{ii} < 0.$$
It follows from Theorem \ref{thm:TM2} that the pair $(A, B)$ is diagonally Riccati stable.  This completes the proof. 

For future use, we note an alternative form of the condition for diagonal Riccati stability given in Theorem \ref{thm:Kra}.  We will use the following simple lemma.

\begin{lemma}\label{lem:DP}
Let $P \in \mathbb{R}^{n\times n}$ be a P-matrix and let $D \succ 0$ be a diagonal matrix in $\mathbb{R}^{n \times n}$.  Then $DPD$ is also a P-matrix.
\end{lemma}
\textbf{Proof:}  Let $x \neq 0$ in $\mathbb{R}^n$ be given.  Then setting $y = Dx$, it follows that there is some index $i$ with $y_i(Py)_i > 0$ as $P$ is a P-matrix.  Thus as $y_i = d_{ii} x_i$ where $d_{ii}$ is the $i$-th diagonal entry of $D$, 
\begin{eqnarray*}
d_{ii}x_i (PDx)_i = x_i(DPDx)_i > 0.
\end{eqnarray*}

\begin{theorem}\label{thm:Kra2}
Let $A, B$ in $\mathbb{R}^{n \times n}$ be given.  The pair $(A, B)$ is diagonally Riccati stable if and only if $-(A\circ S_{11} + B \circ S_{12})$ is a P-matrix for all $S \succeq 0$ in $Sym(2n, \mathbb{R}$ of the form \eqref{eq:HadS} with $\textrm{diag}(S_{11}) = \textrm{diag}(S_{22}) = \mathbf{1}_n$.  
\end{theorem}
\textbf{Proof:} Let $(A, B)$ be diagonally Riccati stable.  Theorem \ref{thm:Kra} implies that $-(A\circ S_{11} + B \circ S_{12})$ is a P-matrix for all $S \succeq 0$ in $Sym(2n, \mathbb{R})$ with $\textrm{diag}(S_{11}) = \textrm{diag}(S_{22}) \gg 0$ so it certainly holds for $S$ where $\textrm{diag}(S_{11}) = \textrm{diag}(S_{22}) = \mathbf{1}_n$.  

Conversely, let $S \succeq 0 $ in $Sym(2n, \mathbb{R})$ with $\textrm{diag}(S_{11}) = \textrm{diag}(S_{22}) \gg 0$ be given.  Then define $D \in \mathbb{R}^{n \times n}$ to be the diagonal matrix with diagonal given by $\textrm{diag}(S_{11})$ and consider $\hat S$ given by $\hat S = TST$ where 
$$T = \left(\begin{array}{c c}
					\sqrt{D^{-1}} & 0 \\
                     0 & \sqrt{D^{-1}}
             \end{array}\right).$$	
It is easy to see that $\hat{S} \succeq 0$ and that $\textrm{diag}(\hat S_{11}) = \textrm{diag}(\hat S_{22}) = \mathbf{1}_n$.  It follows that 
$$-(A \circ \hat{S}_{11} + B \circ \hat S_{12})$$ is a P-matrix.  However,
$$-(A \circ S_{11} + B \circ S_{12}) = -\sqrt{D}(A \circ \hat S_{11} + B \circ \hat S_{12})\sqrt{D}$$
and hence $-(A \circ S_{11} + B \circ S_{12})$ is a P-matrix by Lemma \ref{lem:DP}.  It now follows from Theorem \ref{thm:Kra} that $(A, B)$ is diagonally Riccati stable.  This completes the proof. 

\section{Classes of Diagonally Riccati Stable Matrix Pairs}
In this section, we first note a simple invariance property of diagonal Riccati stability and then show how this result can be used to obtain a wide variety of new classes of diagonally Riccati stable pairs.  

\begin{proposition}
\label{prop:DAD} Let $(A, B)$ be diagonally Riccati stable and let $D$, $E$ be diagonal matrices with diagonal entries $d_{ii}$ and $e_{ii}$ respectively. If $0<e_{ii}^2\leq  d_{ii}^2$, $i=1,\ldots,n$, then $(DAD, DBE)$ is also diagonally Riccati stable. 
\end{proposition}
\textbf{Proof:} Suppose $(A, B)$ is diagonally Riccati stable and let $D$, $E$ be diagonal matrices satisfying the hypotheses of the proposition.  Note, in particular, that both $D$ and $E$ are invertible.  Let $H$ be a non-zero positive semi-definite matrix given by \eqref{eq:H} with $\textrm{diag}(H_{11}) \geq \textrm{diag}(H_{22})$.  Define $\hat{H} = \hat{D} H \hat{D}$ where
$$\hat{D} = \left(\begin{array}{c c}
				D & 0\\
                0 & E
                \end{array}\right).$$
Then it is simple to verify that $\hat{H} \succeq 0$ and is non-zero.  Moreover, it follows from $d_{ii}^2 \geq e_{ii}^2$ that $\textrm{diag}(\hat{H}_{11}) \geq \textrm{diag}(\hat{H}_{22})$.  As $(A, B)$ is diagonally Riccati stable, it follows from Theorem \ref{thm:TM2} that there is some negative diagonal entry of $A\hat{H}_{11} + B \hat{H}_{12}^T$.  This means that the matrix
$$A (DH_{11}D) + B(E H_{12}^T D)$$ has a negative diagonal entry.  Multiplying on the left by $D$ and on the right by $D^{-1}$ doesn't change the sign of any diagonal entries so we conclude that
$$(DAD)H_{11} + (DBE) H_{12}^T$$ has a negative diagonal entry.  Thus $(DAD, DBE)$ is diagonally Riccati stable by Theorem \ref{thm:TM} as claimed.  

We next note that it is also possible to prove the above result directly from the inequality \eqref{eq:Ricc3}; moreover, the following alternative argument also explicitly relates the diagonal matrices solving \eqref{eq:Ricc3} for $(DAD, DBE)$ to the solutions for $(A, B)$. 

\textbf{Alternative proof for Proposition \ref{prop:DAD}.}

Let diagonal positive definite matrices $P$ and $Q$ satisfying \eqref{eq:Ricc3} be given.  Then it follows from the Schur complement \cite{HJ} that 
$$\left(\begin{array}{c c}
	A^TP+PA + Q & PB \\
    B^TP & -Q
    \end{array}\right) \prec 0.
    $$
A simple conjugacy with the matrix 
$$\hat{D} = \left(\begin{array}{c c}
				D & 0 \\
                0 & E
                \end{array}\right)$$
and a little rearrangement using the fact that diagonal matrices commute shows that 
$$\left(\begin{array}{c c}
	(DA^TD)P+P(DAD) + DQD & P(DBE) \\
    (EB^TD)P & -EQE
    \end{array}\right) \prec 0.
    $$
However, as $e_{ii}^2 \leq d_{ii}^2$ for $1\leq i\leq n$, and $Q$ is diagonal, it now follows that
$$\left(\begin{array}{c c}
	(DA^TD)P+P(DAD) + DQD & P(DBE) \\
    (EB^TD)P & -DQD
    \end{array}\right)$$
  $$  \preceq \left(\begin{array}{c c}
	(DA^TD)P+P(DAD) + DQD & P(DBE) \\
    (EB^TD)P & -EQE
    \end{array}\right) \prec 0.
    $$
Hence, by Schur complement again, $P, DQD$, will solve \eqref{eq:Ricc3} for $(DAD, DBE)$. 

The following corollary, which will prove useful in the next subsection is now immediate. 
\begin{corollary}
\label{cor:DAD1} Let $A \in \mathbb{R}^{n \times n}$, $B \in \mathbb{R}^{n \times n}$ be given and let $D$, $E$ be diagonal matrices with $d_{ii} \in \{-1, +1\}$, $e_{ii} \in \{-1, +1\}$ for $1 \leq i \leq n$.  The pair $(A, B)$ is diagonally Riccati stable if and only if $(DAD, DBE)$ is diagonally Riccati stable. 
\end{corollary}

\subsection{Applications}
Proposition \ref{prop:DAD} and Corollary \ref{cor:DAD1} allow us to readily identify classes of diagonally Riccati stable matrix pairs using previous results.  We next provide a (far from exhaustive) list of such classes.  We first introduce some notation necessary for stating our results. 

Let a matrix $C \in \mathbb{R}^{n \times n}$ be given.  We denote by $\hat C$ the matrix with $\hat c_{ii}=c_{ii}$, $\hat c_{ij}=|c_{ij}|$ for $i\neq j$, and we use  $\bar C$ to denote the matrix with entries $\bar c_{ij}=|c_{ij}|$ for  $i,j=1,\ldots,n$.  Note that for any $C$, the matrix $\hat C$ is Metzler while $\bar C$ is nonnegative.  It is known that for $A$ Metzler and $B$ nonnegative, the pair $(A, B)$ is diagonally Riccati stable if and only if $A+B$ is Hurwitz.  We next use Corollary \ref{cor:DAD1} to describe pairs $(A, B)$ for which diagonal Riccati stability is equivalent to the Hurwitz-stability of $\hat A + \bar B$. 

\begin{proposition}
\label{prop:Met1} Let $A \in \mathbb{R}^{n \times n}$ be Metzler and $B = e_k b^T$ for some $b \in \mathbb{R}^n$ and some $k \in \{1, \ldots, n\}$.  Then $(A,B)$ is diagonally Riccati stable if and only if $\hat A + \bar B$ is Hurwitz.
\end{proposition}
\textbf{Proof:} Take $D = I$ and define $E$ by setting $e_{ii} = \textrm{sign}(b_i)$ for $1 \leq i \leq n$.  It is simple to see that $\bar B = DBE$ and $\hat A = DAD$.  The result is now a simple application of Corollary \ref{cor:DAD1}. 

Our next result concerns matrices $B$ of the same form as in the previous proposition with sign-symmetric tridiagonal matrices $A$ which need not be Metzler. 

\begin{proposition}\label{prop:Met2}
Let $A \in \mathbb{R}^{n \times n}$ be of the form
\begin{equation}\label{eq:A1}
A = \begin{pmatrix}
a_{1} & u_{1} & 0 & \cdots & 0 & 0 \\
l_{1} & a_{2}&  u_{2}  & \cdots & 0 &0 \\    
0 & l_{2} & a_{3} & \cdots & 0 & 0\\ 
\vdots & \vdots & \vdots & \ddots & \vdots & \vdots \\
0 & 0 & 0 & \cdots  & a_{n-1} & u_{n-1}\\
0 & 0 & 0 & \cdots  & l_{n-1} &  a_{n}\\
\end{pmatrix},               
\end{equation}
where 
$$
l_i u_i\geq 0, \quad i=1,\ldots,n-1. 
$$  Let $B = e_k b^T$ for some $b \in \mathbb{R}^n$ and some $k \in \{1, \ldots, n\}$.  Then $(A, B)$ is diagonally Riccati stable if and only if $\hat A + \bar B$ is Hurwitz.  
\end{proposition}
\textbf{Proof:} We again use Corollary \ref{cor:DAD1}.  We first define the matrix $D$ by setting $d_{11} = 1$ and $$d_{ii} = \frac{\textrm{sign}(l_{i-1})}{d_{i-1, i-1}}$$ for $2 \leq i \leq n$.  It is then easy to see that $\hat A = DAD$.  Next we define the diagonal matrix $E$ by setting $e_{ii} = \textrm{sign}(d_{ii}b_i)$ for $1 \leq i \leq n$.  Then $\bar B = DBE$ and it follows again from Corollary \ref{cor:DAD1} that $(A, B)$ is diagonally Riccati stable if and only if $\hat A + \bar B$ is Hurwitz. 

A similar result also holds for matrices $A$ of the form:
\begin{equation}\label{eq:A2}
A = \begin{pmatrix}
a_{1} & 0 & 0 & \cdots & 0 & 0 \\
0 & a_{2}&  0  & \cdots & 0 &0 \\    
0 & 0 & a_{3} & \cdots & 0 & 0\\ 
\vdots & \vdots & \vdots & \ddots & \vdots & \vdots \\
0 & 0 & 0 & \cdots  & a_{n-1} & 0\\
c_{1} & c_{2} & c_{3} & \cdots  & c_{n-1} &  a_{n}\\
\end{pmatrix}.                       
\end{equation}
Here, consider $D$ given by $d_{nn} = 1$ and $d_{ii} = d_{i+1, i+1}\textrm{sign}(c_i)$ for $1 \leq i \leq n-1$.  The $DAD = \hat A$ and if we define $E$ as in the proof of Proposition \ref{prop:Met2}, we will have $\bar B = DBE$.  Corollary \ref{cor:DAD1} implies the following result.
\begin{proposition}
\label{prop:Met3} Consider $A$ in $\mathbb{R}^{n \times n}$ given by \eqref{eq:A2} and $B = be_k^T$ for some $b \in \mathbb{R}^n$ and some $k \in \{1, \ldots, n\}$.  The pair $(A, B)$ is diagonally Riccati stable if and only if $\hat A + \bar B$ is Hurwitz.  
\end{proposition}

Finally for this section, we note that analogous results to those given above can be obtained for the case where $B$ is of the form
\begin{equation}\label{eq:B2}
B =\begin{pmatrix}
0 & b_{1} & 0 & \cdots & 0 & 0 \\
0 & 0&  b_{2}  & \cdots & 0 &0 \\    
0 & 0 & 0 & \cdots & 0 & 0\\ 
\vdots & \vdots & \vdots & \ddots & \vdots & \vdots \\
0 & 0 & 0 & \cdots  & 0 & b_{n-1}\\
0 & 0 & 0 & \cdots  & 0 &  0\\
\end{pmatrix}.
\end{equation}
If we choose a diagonal $D$ (with $d_{ii}= \pm 1$ for $1 \leq i \leq n$) such that $DAD = \hat A$, then we can ensure that $DBE = \bar B$ by taking $E$ to be the diagonal matrix with $e_{ii} = \textrm{sign}(d_{i-1, i-1} b_{i-1})$ for $2 \leq i \leq n$ and $e_{11} = 1$.  Using the simple observation, the following result can be established identically to Propositions \ref{prop:Met1}, \ref{prop:Met2}, \ref{prop:Met3}. 

\begin{proposition}\label{prop:Met4}
Let $A$ in $\mathbb{R}^{n \times n}$ be Metzler or in one of the forms \eqref{eq:A1}, \eqref{eq:A2}.  If $B$ is of the form \eqref{eq:B2}, then the pair $(A, B)$ is diagonally Riccati stable if and only if $\hat A + \bar B$ is Hurwitz.
\end{proposition}

\section{Classes of Diagonally Riccati Stable $3 \times 3$ Matrix Pairs}
All of the classes of matrix pairs analysed in the last section have the property that there exist diagonal matrices $D,E$ satisfying the hypotheses of Corollary \ref{cor:DAD1} such that $DAD$ is Metzler and $DBE$ is nonnegative. In this section, we consider classes of matrix pairs for which this is not possible and provide necessary and sufficient conditions for diagonal Riccati stability for these.  Specifically, we derive such results for two classes of matrix pairs in $\mathbb{R}^{3 \times 3}$ which arise in the study of indirect control systems. 

Consider the matrices
\begin{equation}\label{eq:3AB1}
A= \begin{pmatrix}
a_1 & 0 & 0  \\
c_1 & a_2&  0 \\    
0 & c_2 & a_3 \\
\end{pmatrix}, \quad
B= \begin{pmatrix}
0 & 0 & b_1  \\
0 & 0&  b_2 \\    
0 & 0 & 0 \\
\end{pmatrix}.                      
\end{equation}

We will need the following technical lemma, which can be readily verified by the method of Lagrange multipliers.

\begin{lemma}\label{lem:Lag}
Let $C, D \in \mathbb{R} \backslash \{0\}$ be given.  Consider $$F(x, y, z) = Cx + D yz.$$ 
The value of $|F|$ subject to $x, y, z \in [-1, 1]$ and 
\begin{equation}
\label{eq:LemLag1} 1 - (x^2 + y^2 + z^2) + 2xyz \geq 0
\end{equation}
is bounded by $\max \{|C|, |D+C|\}$. 
\end{lemma}
\textbf{Proof:}  If $x \in\{-1, 1\}$, then it is not difficult to see that $y, z$ are both $0$ so that $|F| = |C|$.  On the other hand, if either of $|y|, |z|$ is 1, then the other two variables must be zero and $|F| = 0$.  From the form of $F$ it is easy to see that the extrema of $F$ (and hence of $|F|$) must be at points where the inequality in \eqref{eq:LemLag1} is an equality.  Thus we can use the method of Lagrange multipliers.  The Lagrangian is given by 
$$\mathcal{L} = Cx + D yz + \lambda (1 - (x^2 + y^2 + z^2) + 2xyz).$$
If we consider the conditions given by $\frac{\partial{\mathcal{L}}}{\partial y} = 0$, $\frac{\partial{\mathcal{L}}}{\partial z} = 0$, we find that 
$$Dz - 2\lambda y + 2 \lambda xz = 0,$$
$$Dy - 2\lambda z + 2 \lambda xy = 0.$$
Moreover, the condition given by the partial derivative with respect to $x$ shows that $\lambda \neq 0$.   Multiplying the first equation by $y$, the second by $z$ and subtracting, we see that $y^2 = z^2$ so that at any extremum, either $y =z$ or $y = -z$.  Using the equality in \eqref{eq:LemLag1} it follows that either $x = 1- 2y^2$, $y= -z$ or $x = 2y^2 - 1$, $y = z$.  In either of these cases, the absolute value of $F$ is given by 
$$|(2C + D)y^2 - C|.$$  As $|y| \leq 1$, it follows that $|F|$ is bounded above by $\max\{|C+D|, |C|\}$ as claimed. 

\begin{theorem}
\label{thm:5}
Let $A$ and $B$ be given by \eqref{eq:3AB1}.  The pair $(A, B)$ is diagonally Riccati stable if and only if:
\begin{itemize}
\item[(i)] $a_i < 0$ for $1 \leq i \leq 3$;
\item[(ii)] $a_2a_3 > |b_2c_2|$;
\item[(iii)] $|a_1a_2a_3| > |c_2(b_1c_1 - a_1b_2)|$.
\end{itemize}
\end{theorem}

{\bf Proof.}  We will make use of Theorem \ref{thm:Kra2} so let $S \succeq 0$ in $Sym(6, \mathbb{R})$ be given by
\begin{equation}\label{eq:SUV}
S = \left(\begin{array}{c c}
		U & W\\
        W^T & V
        \end{array}\right)
\end{equation}
with $\textrm{diag}(U) = \textrm{diag}(V) = \mathbf{1}_3$.  It follows easily from the fact that all $2 \times 2$ principal minors of $S$ are nonnegative that for all $i, j$, $|u_{ij}| \leq 1$, $|w_{ij}| \leq 1$.  Moreover the matrix $A\circ U + B \circ W$ is given by
\begin{equation}\label{eq:ABP1}
T :=\left(\begin{array}{c c c}
		a_1 & 0 & w_{13}b_1 \\
        u_{21}c_1 & a_2 & w_{23}b_2\\
        0 & u_{32}c_2 & a_3
        \end{array}\right).
\end{equation}
If $(A, B)$ is diagonally Riccati stable, Theorem \ref{thm:Kra2} implies that $-T$ must be a P-matrix for the case where the matrix $S$ has $U = V = W = \mathbf{1}_3 \mathbf{1}_3^T$ and the case where $U = V = \mathbf{1}_3 \mathbf{1}_3^T$ and $W = - U$.  Checking the minors of order 1, 2 and 3 in these cases yields conditions (i), (ii) and (iii).  

For the converse, assume that (i), (ii), (iii) hold and let $S \succeq 0$ of the form \eqref{eq:SUV} be given.  Then as noted above 
\begin{equation}
\label{eq:uwLag1} |u_{ij}| \leq 1, \; |w_{ij}| \leq 1 \forall i, j.
\end{equation}
Moreover, if we consider the $3 \times 3$ submatrix of $S$ formed from rows and columns with indices in $\{1, 2, 6\}$ we see that 
\begin{equation}\label{eq:uwLag2}
1 - (u_{12}^2 + w_{23}^2 + w_{13}^2) + 2u_{12}w_{13}w_{23} \geq 0.
\end{equation}

It is enough to show that for the matrix $T$ in \eqref{eq:ABP1}, $-T$ is a P-matrix provided \eqref{eq:uwLag1}, \eqref{eq:uwLag2} are satisfied.  We need to check all principal minors of $-T$ are positive.  Condition (i) implies that this is true for the minors of size 1.  The 3 principal minors of size 2 are given by 
$$a_1a_2, \, a_1a_3, \, a_2a_3 - u_{32}w_{23}b_2c_2.$$
Conditions (i) and (ii) together imply that all 3 of these are positive.  Finally, the principal minor of $T$ size 3 (determinant) is given by
\begin{eqnarray}\nonumber
& & a_1a_2a_3 - w_{23}u_{32}a_1b_2c_2 + w_{13}u_{21}u_{32}b_1c_1c_2\\
&=& a_1a_2a_3 +u_{32}c_2(w_{13}u_{21}b_1c_1-w_{23}a_1b_2).
\end{eqnarray}
If either of $b_1c_1$ or $a_1b_2$ is zero, then it follows readily from (iii) combined with (i) that $\textrm{det}(T) < 0$ and that $-T$ is a P-matrix.  On the other hand, if both $b_1c_1$ and $a_1b_2$ are non-zero, it follows from Lemma \ref{lem:Lag} that the extremal values of $|w_{13}u_{21}b_1c_1-w_{23}a_1b_2|$ subject to \eqref{eq:uwLag1}, \eqref{eq:uwLag2} must be one of $|a_1b_2|$ or $|a_1b_2 - b_1c_1|$.  In the former case, it follows from (ii) and (i) that $\textrm{det}(T) < 0$, while in the latter, the same conclusion follows from combining (i) and (iii).  This completes the proof.

For our final result, we consider matrices in $\mathbb{R}^{3 \times 3}$ of the form:
\begin{equation}\label{eq:3AB2}
A = \begin{pmatrix} a_1 & 0 & 0 \\0 & a_2 & 0  \\ c_1 & c_2 & a_3  
\end{pmatrix},             \qquad
B = \begin{pmatrix} 0 & 0 & b_1\\0 & 0 & b_2\\ 0 & 0 & 0 \end{pmatrix}.
\end{equation}

\begin{theorem}
\label{thm:6}
Let $A$ and $B$ be given by \eqref{eq:3AB2}.  The pair $(A, B)$ is diagonally Riccati stable if and only if:
\begin{itemize}
\item[(i)] $a_i < 0$ for $1 \leq i \leq 3$;
\item[(ii)] $a_1a_3 > |c_1b_1|$, $a_2a_3 > |b_2c_2|$;
\item[(iii)] $|a_1a_2a_3| > |a_1b_2c_2+b_1c_1a_2|$.
\end{itemize}
\end{theorem}

{\bf Proof.}  
The argument follows a very similar path to that of Theorem \ref{thm:5} and again relies on Theorem \ref{thm:Kra2}.  Let $S \succeq 0$ in $Sym(6, \mathbb{R})$ be given by \eqref{eq:SUV}
with $\textrm{diag}(U) = \textrm{diag}(V) = \mathbf{1}_3$.  Once again, it follows that for all $i, j$, $|u_{ij}| \leq 1$, $|w_{ij}| \leq 1$.  

In this instance, the matrix $A\circ U + B \circ W$ is given by
\begin{equation}\label{eq:ABP2}
T :=\left(\begin{array}{c c c}
		a_1 & 0 & w_{13}b_1 \\
        0 & a_2 & w_{23}b_2\\
        u_{31}c_1 & u_{32}c_2 & a_3
        \end{array}\right).
\end{equation}
If $(A, B)$ is diagonally Riccati stable, (i), (ii), (iii) again follow from Theorem \ref{thm:Kra2}: considering separately the cases corresponding the the matrix $S$ where $U = V = W = \mathbf{1}_3 \mathbf{1}_3^T$ and the case where $U = V = \mathbf{1}_3 \mathbf{1}_3^T$ and $W = - U$.  

For the converse, we show that $-T$ is a P-matrix where the matrix $T$ is given by \eqref{eq:ABP2} and $|u_{ij}| \leq 1$, $|w_{ij}| \leq 1$ for all $i, j$.  Conditions (i) and (ii) imply that all principal minors of $T$ of size 1, 2 are negative, positive respectively.  Finally, the principal minor of size 3 (determinant) is given by
\begin{eqnarray*}
& & a_1(a_2a_3 - w_{23}u_{32}b_2c_2) - w_{13}u_{31}b_1c_1a_2\\
&=& a_1a_2a_3 - (u_{23}w_{23} a_1b_2c_2 + w_{13}u_{13}b_1 c_1 a_2).  
\end{eqnarray*}
If we consider separately the submatrices of $S$ formed from rows/columns with indices in $\{1, 3, 6\}$, $\{2, 3, 6\}$ we see that:
\begin{equation}\label{eq:Lag3}
1-w_{33}^2 - (u_{13}^2 + w_{13}^2) + 2 w_{33}u_{13}w_{13} \geq 0;
\end{equation}
\begin{equation}\label{eq:Lag4}
1-w_{33}^2 - (u_{23}^2 + w_{23}^2) + 2 w_{33}u_{23}w_{23} \geq 0.
\end{equation}
As in the proof of Theorem \ref{thm:5} we shall examine the extrema of  
$$F(u_{13}, w_{13}, u_{23}, w_{23}) := u_{23}w_{23} a_1b_2c_2 + w_{13}u_{13}b_1 c_1 a_2$$ subject to \eqref{eq:Lag3}, \eqref{eq:Lag4} and the additional constraints that all entries of $S$ are in $[-1, 1]$.  

We first note that if $a_1b_2c_2$, $b_1c_1a_2$ have the same sign, then clearly
$$|F(u_{13}, w_{13}, u_{23}, w_{23})| \leq  |a_1b_2c_2+b_1c_1a_2|.$$
It follows from (i), (ii) and (iii) that in this case $\textrm{det}(T) < 0$ so $-T$ is a P-matrix.  

Now suppose that $a_1b_2c_2$, $b_1c_1a_2$ have opposite signs.  For ease of notation, we will write $C =  a_1b_2c_2$, $D = b_1c_1a_2$, $\alpha = w_{33}$, $x = u_{23}$, $y = w_{23}$, $z = u_{13}$, $w = w_{13}$.  It is not difficult to see that at an extremum of $F$, the inequalities in \eqref{eq:Lag3}, \eqref{eq:Lag4} must be tight as otherwise we could increase or decrease the absolute value of $F$ by suitably altering one of the pairs $(x, y)$, $(z, w)$ without changing the other.    

So we wish to find the extreme values of $|F(x, y, z, w)| = |Cxy + Dzw|$ subject to $|x| \leq 1$, $|y| \leq 1$, $|z| \leq 1$, $|w| \leq 1$ and 
\begin{equation}\label{eq:Lagxy}
(1-\alpha^2) - (x^2+y^2) + 2\alpha xy = 0;
\end{equation}
\begin{equation}\label{eq:Lagzw}
(1-\alpha^2) - (z^2+w^2) + 2\alpha zw = 0.
\end{equation}
As $C$ and $D$ have opposite signs, if $xy$, $zw$ have the same sign then $|F| \leq \max\{|C|, |D|\}$.  If we use Lagrange multipliers to find the extrema of $Cxy$ subject to \eqref{eq:Lagxy} it is not difficult to see that either $x = y$ or $x = -y$.  Using \eqref{eq:Lagxy}, the corresponding values of $Cxy$ are given by $C\frac{1+\alpha}{2}$, $C\frac{1-\alpha}{2}$.  The same conclusion holds for the extrema of $Dwz$ subject to \eqref{eq:Lagzw} and the corresponding values of $Dwz$ are $D\frac{1+\alpha}{2}$, $D\frac{1-\alpha}{2}$.  The only case left to consider is when $xy$ and $zw$ have opposite signs.  It is now easy to see that in this case we also have $$|F| \leq \max\{|C|, |D|\} = \max\{|a_1b_2c_2|, |b_1 c_1 a_2| \}.$$  This together with (i), (ii) implies that $\textrm{det}(T) < 0$ and hence $-T$ is a P-matrix in this case also.  This completes the proof.

\section{Conclusions}
We have presented an extension of Kraiijevanger's condition for Lyapunov diagonal stability to Riccati stability and time-delay systems.  We have also shown how diagonal Riccati stability of a time-delay system is invariant under certain transformations on the defining matrix pair.  These results have then been used to provide simple conditions for diagonal stability for certain classes of time-delay systems.  In future work, the authors wish to use the work presented here to extend results such as those found in \cite{Wimmer, Arcak} to time-delay systems.

\section*{Acknowledgments}
This work was partially supported by the Russian Foundation for Basic Research, grant no. 16-01-00587 and by Science Foundation Ireland grant 13/RC/2094 and co-funded under the European Regional Development Fund through the Southern \& Eastern Regional Operational Programme  to Lero - the Irish Software Research Centre (www.lero.ie).

\end{document}